\documentclass{amsart}
\usepackage{amssymb}
\usepackage{stmaryrd}
\usepackage{booktabs}

\newcommand{\ZZ}{\mathbf{Z}}
\newcommand{\QQ}{\mathbf{Q}}
\newcommand{\ep}{\varepsilon}
\newcommand{\divides}{\mathbin{|}}
\newcommand{\ps}[1]{\llbracket #1 \rrbracket}

\DeclareMathOperator{\Gal}{Gal}
\DeclareMathOperator{\Cyc}{Cyc}

\title{Irregular primes to 163 million}
\author{J.~P.~Buhler}
\author{D.~Harvey}

\begin{document}

\begin{abstract}
We compute all irregular primes less than $163\,577\,356$.  For all of these primes we verify that the Kummer--Vandiver conjecture holds and that the $\lambda$-invariant is equal to the index of irregularity.
\end{abstract}

\maketitle

\section{Introduction}

Bernoulli numbers are deeply intertwined with the
the arithmetic of cyclotomic fields.   To explain a basic
case, consider the cyclotomic field $K = \QQ(\zeta_p)$ generated by 
a $p$-th root of unity $\zeta_p = e^{2\pi i/p}$, $p$~prime. Let 
$G = \Gal(K/\QQ) \simeq (\ZZ/p\ZZ)^\times$ be its Galois group, and let 
$\omega : G \to \ZZ_p^\times$ be the unique $p$-adic character for which 
\[
\omega(\sigma) \equiv r \bmod p
\]
where $\sigma(\zeta_p) = \zeta_p^r$.  Then the $p$-Sylow subgroup $A$ of
the class group decomposes into components
\[
A = \bigoplus_{0 < k < p-1} A_k
\]
where $A_k$ is the $\ZZ_p[G]$-submodule on which $\sigma$ acts
by multiplication by $\omega^k(\sigma)$.  (For details, see \cite{Washington}.)

It is conjectured that $A_k$ is trivial if $k$ is even; this is equivalent to the
Vandiver (or Kummer--Vandiver) conjecture that the class number $h^+_p$ of the 
maximal real subfield $\QQ(\zeta_p+\zeta_p^{-1})$ is not divisible by~$p$.  
By combining results of Herbrand and Ribet 
we know that $A_{p-k}$ is nontrivial if and only if the numerator of the $k$-th
Bernoulli number is divisible by~$p$, and the Main Conjecture, proved by Mazur
and Wiles, shows that the order of $A_{p-k}$ is equal to the power of~$p$
that divides~$B_k$.
Thus, the first step in finding the $p$-part of the class group of~$K$
is to study the numerators of Bernoulli numbers.

If $k$ is an even integer with $0 < k < p-1$, then $k$ is said to be
an \emph{irregular index} for~$p$ if $p$ divides the numerator of the $k$-th 
Bernoulli number~$B_k$.  
The \emph{index of irregularity} of~$p$, denoted $i_p$, is the number 
of such irregular indices~$k$, and $p$ is \emph{irregular} if $i_p$ 
is positive.  If $k$ is irregular for~$p$ then $(p,k)$ is said to 
be an \emph{irregular pair}.  Once the irregular pairs for a given $p$ have been 
found it is possible to make further calculations verifying the Kummer--Vandiver 
conjecture, and to calculate the corresponding $\lambda$-invariants 
in Iwasawa theory; these invariants $\lambda_p$ measure the rate of growth
of the $p$-part of the class group of $p$-power cyclotomic fields.

This paper describes a calculation of all irregular pairs
for $p$ less than $39 \cdot 2^{22} = 163\,577\,356$. 
The computation consumed more than twenty years of CPU time, and the quick summary
of the results is that nothing surprising happened.
The Kummer--Vandiver conjecture
is true for all of these primes, and the $\lambda$-invariant
is always as small as it can be, i.e., $\lambda_p = i_p$.  This implies
that the class group of $\QQ(\zeta_{p^n})$ is, for all of these~$p$, isomorphic to
$(\ZZ/p^n \ZZ)^{i_p}$.  We did find one new prime with $i_p = 7$, 
namely $p = 32\,012\,327$. At the moment that we write this, the data supporting
these calculations can be found at \cite{URL}.

Irregular indices have been computed many times over the 
last 160 years, starting with Kummer's hand calculations for $p \le 163$,
which he used to verify Fermat's Last Theorem and that $h^+_p$ is prime to~$p$.
This was extended to $p \le 619$ early in the twentieth century by Vandiver and
his colleagues, using desk calculators (and graduate students).  Vandiver and
the Lehmers, making one of the first number-theoretic calculations on electronic
computers, extended this to $p \le 4001$ by 1956. For a discussion of these
early results, see \cite{Corry}.  A series of further results by the Lehmers,
Iwasawa, Sims, Kobelev, Johnson, Nicol, Pollack, Selfridge, Tanner, and Wagstaff 
extended these to $p < 150\,000$ by 1987 (see \cite{John}, \cite{Wag}, or \cite{WagTan}
for a summary and further references).
In the early 1990's,
Buhler, Crandall, Ernvall,  Mets{\"a}nkyl{\"a}, Sompolski, and Shokrollahi used
``asymptotically fast'' algorithms to go much further, and by 1999 the
computations were extended to all primes less than 12 million 
\cite{1-million},
\cite{4-million},
\cite{12-million},
\cite{EMa},
\cite{EMb}.

This paper describes our extension of these results.
The two major reasons that we were able to extend the upper bound on $p$ by a factor
of almost 14 were: 
(a) the use of
two large clusters at the Texas Advanced Computing Center at the 
University of Texas at Austin, and
(b) new algorithms and implementations of them, due to the second author, for 
polynomial arithmetic over finite fields.

The paper is organized as follows: we first describe
the algorithms for calculating $B_k \bmod p$ for $0 < k < p-1$,
algorithms for performing the needed polynomial arithmetic, and
the algorithms for the Vandiver and cyclotomic invariant calculations.
Then we summarize the results, including a description of the hardware
that we used and a discussion of the steps we took to ensure correctness.

\section{Bernoulli numbers modulo $p$}

Bernoulli numbers are defined by a formal power series inversion
\[
\frac{t}{e^t-1} = \sum_{k \ge 0} B_k \frac{t^k}{k!}.
\]
We have $B_1 = 1/2$, and all other odd Bernoulli numbers $B_3, B_5, \ldots$ are zero.
The first, and by far and away the most time-consuming phase of our calculations is
the computation of $B_k \bmod p$, $0 \le k < p-1$, $k$ even, for each~$p$ in turn.

We used two different algorithms, the \emph{Voronoi identity method}
and the \emph{power series method}. Both depend crucially on
asymptotically fast polynomial arithmetic over $\ZZ/p\ZZ$. The
theoretical complexity is $O(p \log^{2+\ep} p)$ for both methods,
but the implied constants and memory footprints are different --- the algorithm
using the Voronoi identity seems to be faster but requires more memory. We
first describe the two methods, and then discuss the algorithms
used for the underlying polynomial arithmetic.

The method based on the Voronoi identity is essentially 
identical to the `root finding method' of \cite{12-million}; we used it 
for all $p < 21 \cdot 2^{22} = 88\,080\,384$. Let $g$ be a generator of 
$(\ZZ/p\ZZ)^\times$. Then
for $2 \leq k < p-1$, $k$ even, we have

\begin{equation}
\label{eq:voronoi}
 B_{k} \equiv \frac{2k}{1 - g^{k}} \sum_{i=0}^{(p-3)/2} g^{(k-1)i} h(g^i) \pmod p,
\end{equation}
where
 \[ h(x) = \frac {(x \bmod p) - g(x/g \bmod p)}p + \frac{g-1}2. \]
This is a version of Voronoi's congruence with the terms 
arranged in multiplicative order; see \cite{Slav}, \cite{Vnd} (or \cite{Dil} for a Bernoulli number 
bibliography, with an updated online version at \cite{DilURL}).
The sum on the right
hand side may be regarded as the Fourier transform of 
$i \mapsto h(g^i)/g^i$ with respect to the roots of unity
$\{g^{2k}\}_{k=0}^{(p-3)/2}$. Using Bluestein's FFT algorithm
\cite{bluestein}, we convert the problem of simultaneously computing
all $B_k$ to the problem of multiplying two polynomials of length
$(p-1)/2$ over $\ZZ/p\ZZ$, plus $O(p)$ pre-processing and post-processing
operations in $\ZZ/p\ZZ$.

The power series method is described in \cite{1-million}. We used
the variant based on the identity
 \[ \sum_{n=0}^\infty \frac{2^n B_{2n}}{(2n)!} x^{2n} = \frac{A_0(x^8)
 + A_2(x^8)x^2 + A_4(x^8)x^4 + A_6(x^8)x^6}{D(x^8)}, \]
where $D$ and the $A_i$ are power series whose coefficients are
given by simple recurrence formulae \cite[p.~719]{1-million}.
This reduces to a power series inversion of length $p/8 + O(1)$, 
and four series multiplications of length $p/8 + O(1)$,
over $\ZZ/p\ZZ$. This method uses less memory than the Voronoi
identity method, and we used it for all primes $21 \cdot 2^{22} <
p < 39 \cdot 2^{22}$.

\section{Polynomial and power series arithmetic}

Most of the algorithms for polynomial and power series arithmetic described here 
are implemented in the zn\_poly library of the second author 
\cite{zn-poly}.

In what follows, we consider arithmetic in $R[x]$, including
arithmetic on truncated power series in $R \ps x$, where $R$ is a
coefficient ring in which $2$ is not a zero-divisor. For the algorithms
described above, we could take $R = \ZZ/p\ZZ$, where $p$ is the
prime of interest. In our implementation, we found it profitable
to handle two nearby primes simultaneously, taking instead $R =
\ZZ/p_1 p_2\ZZ$, and using the Chinese Remainder Theorem. This was
possible since the computations were run on 64-bit hardware, and
the largest prime considered was less than 32 bits long.

Our large polynomial multiplication routine involves a combination
of algorithms. The outermost layer reduces the multiplication of
polynomials $f, g \in R[x]$ of length $n$ to $O(n^{1/2})$ negacyclic
convolutions of length $O(n^{1/2})$, using an adaptation to polynomials
of the Sch\"onhage--Strassen integer multiplication algorithm
\cite{schonhage-strassen} suggested by Sch\"onhage \cite{schonhage}.
This proceeds by splitting the inputs into segments of length $M/2$
and embedding the multiplication problem into $R[y,z]/(y^M+1, z^K-1)$
via the inverse of the substitution $y \mapsto x$, $z \mapsto
x^{M/2}$, where $M = 2^{\lceil(1 + \log_2 n)/2\rceil} = O(n^{1/2})$
and $K = 2^{\lceil \log_2(4n/M)\rceil} = O(n^{1/2})$. These conditions
ensure that $MK \geq 4n$ (the product `fits' into the bivariate
quotient ring) and that $K \divides 2M$. The bivariate multiplication
may then be effected by using FFTs over the ring $R[y]/(y^M + 1)$
with respect to the primitive $K$-th root of unity $y^{2M/K}$. These FFTs involve only additions and subtractions in $R$,
plus some book-keeping to keep track of powers of $y$. To improve
memory locality and smoothness of the running time with respect to
$n$, we use a cache-friendly version \cite{cache-trunc-fft} of van
der Hoeven's truncated Fourier transform (TFT) method \cite{tft1, tft2}.

For the negacyclic convolutions (recursive multiplications in $R[y]/(y^M+1)$) we use Nussbaumer's algorithm
\cite{nussbaumer} for sufficiently large $M$; this follows the same
general plan as Sch\"onhage's algorithm sketched above, but uses a
slightly different substitution that can reduce the
size of the Fourier coefficients by a factor of two in some cases.
For sufficiently small $M$, we switch to a variant of the Kronecker
substitution method \cite{multipoint}, which packs the input
polynomials into integers and multiplies the resulting integers.

Our power series inversion algorithm is based on a Newton iteration.
We followed the approach suggested in 
\cite[Algorithm MP-inv, p.~421]{middle-product}, 
which takes advantage of the middle product
to improve the running time constant. A minor improvement in our
algorithm is that we reuse the Fourier transform of $\alpha$ in
lines 3 and 4 of \cite{middle-product}, reducing the theoretical
cost by a further 1/6. To transport the smoothness benefits of the
truncated FFT from the ordinary product to the middle product, it
was necessary to implement a \emph{transposed} TFT and inverse TFT
over $R[y]/(y^M + 1)$, making use of the ``transposition principle'' \cite{tellegen}. Unlike the usual (non-truncated) FFT, the
matrices for the TFT and inverse TFT are not symmetric, so genuinely
new code was required.

It is worth pointing out that our polynomial arithmetic routines
are based entirely on integer arithmetic; we did not use floating
point arithmetic at all.

Memory constraints became formidable for those $p$ near the top of
the range that we handled using the Voronoi identity method. Just
to store two polynomials of length $p/2$ and their product, where
$p \approx 21 \cdot 2^{22}$, requires 1.4 GB of memory, assuming
a 64-bit word for each coefficient. Each thread had only 2 GB RAM
available (see description of the hardware below), leaving very little auxiliary working
storage, certainly much less than is required to store the FFT
coefficients. To overcome this, we implemented a discrete Fourier
transform variant that uses a naive quadratic-time DFT algorithm
to handle the first few layers of the transform, and then the
standard $O(n \log n)$ algorithm for the remaining layers. The resulting time-space tradeoff was acceptable up to some limit,
but became infeasible as we approached the 2 GB barrier, prompting
us to switch to the less memory-intensive power series method.

\section{Kummer--Vandiver and cyclotomic invariants}

We quickly summarize some well-known facts about cyclotomic
class groups, using notation established above; for details see \cite{Washington}.
If $p$ is prime and $k$ is even, $0 \le k < p-1$, then 
$A_{p-k}$ is nontrivial if and only if $k$ is irregular 
for~$p$.  Moreover, a ``reflection theorem'' implies that if
$A_k$ is nontrivial then $A_{p-k}$ is nontrivial, i.e., $k$ is
irregular for~$p$.
The component $A_k$ has the same order as the $\omega^k$-component of the $p$-Sylow subgroup of $U_K/\Cyc_K$, 
where $U_K = O_K^\times$ is the group of units inside the ring of integers of~$K$, and $\Cyc_K$ denotes the 
group of cyclotomic units.
The latter can be described explicitly, and 
Vandiver's conjecture for~$p$, which says that $A_k$ is trivial for
all even~$k$, can be proved by showing that the cyclotomic unit
$u_k$ that generates this component
does not have a $p$-th root lying in~$K$.

This can be established by finding a prime ideal $Q$ of $K$ that
has residue class degree 1, for which $u_k$ is not a $p$-th
power modulo~$Q$.  This reduces (again, for details see \cite{Washington})
to showing that it is possible to find a rational prime~$q$ with $q \equiv 1 \bmod p$
such that $V_{p,k}^{(q-1)/p} \not\equiv 1 \bmod q$, where
 \[ V_{p, k} = \prod_{c=1}^{(p-1)/2} (z^c - z^{-c})^{c^{p-1-k}}
 \pmod q. \]
Here $z$ is a $p$-th root of unity modulo $q$, and $Q$ is one of the primes
lying over~$q$.  Thus, Vandiver's
conjecture for~$p$ is proved if the above condition is checked for
each irregular index $k$ for~$p$.  We performed this computation
for every irregular pair in our table, and in each case the above
condition held, and in all cases the smallest prime $q \equiv 1 \bmod p$ 
worked.  

To verify that the $\lambda$-invariant $\lambda_p$ is equal to $i_p$ it suffices
to check that certain congruences between Bernoulli numbers do not hold.
These reduce to elementary congruences \cite{EMa}.
To describe these, let
\[
S(e) = \sum_{a=1}^{(p-1)/2} a^e.
\]
For each irregular index $k$ for a given prime $p$ we check that
\begin{align*}
2^k & \not\equiv 1 \bmod p, \\
S(k-1) & \not\equiv S(p+k-2) \bmod p^2, \\
k S(k-1) & \not\equiv (k-1)S(p+k-2) \bmod p^2.
\end{align*}
These congruences imply that $B_k \not\equiv 0 \bmod p^2$,
and if they hold for all irregular indices $k$ for~$p$ then $\lambda_p = i_p$.
In the case that $2^k \equiv 1 \bmod p$, this test can be replaced
by another one; for details on this see \cite{EMa} (though note that
our notation is different).

\section{Results}

There are $N = 9\,163\,831$ primes up to $39\cdot 2^{22} = 163\,577\,356$.
The indices of irregularity $i_p$ ranged up to 7, and the number
of primes $N_i$ of index~$i$ are tabulated below.

\begin{table}[h]
\begin{tabular}{crllr}
\toprule
$i$ &   $N_i$      & $N_i/N$      & $p_i$       & $N \cdot p_i$  \\
\midrule
0   & 5\,559\,267  & 0.6067       & 0.6065      & 5\,558\,144    \\
1   & 2\,779\,293  & 0.3033       & 0.3032      & 2\,779\,072    \\
2   &    694\,218  & 0.0758       & 0.0758      &    694\,768    \\
3   &    115\,060  & 0.01256      & 0.01263     &    115\,794    \\
4   &     14\,425  & 0.00157      & 0.00158     &     14\,474    \\
5   &      1\,451  & 0.000158     & 0.000158    &      1\,447    \\
6   &         112  & 0.000012     & 0.000013    &         120    \\
7   &           5  & 0.00000055   & 0.00000094  &           8    \\
\bottomrule
\end{tabular}
\caption{Irregular index statistics for $p < 39 \cdot 2^{22}$}
\end{table}

Here $p_i = e^{-1/2}/(2^i i!)$ is the probability that a Poisson process with mean $1/2$ has $i$ successes.
This is motivated by the heuristic (Lehmer, Siegel) that if the Bernoulli
numerators are random modulo $p$, then $i_p$ is the number of successes in $p/2$ trials
each having probability~$1/p$ of success.

The five primes of index 7 are 3\,238\,481, 5\,216\,111, 5\,620\,861, 9\,208\,289 and 32\,012\,327.

The Vandiver conjecture was found to be true for all primes tested.
The $\lambda$-invariant was equal to $i_p$ for all primes; the alternate
test mentioned above was required for 3 primes.

For a discussion of heuristics that might apply to Vandiver's conjecture
and to the value of $\lambda_p$ see \cite[p.~159]{Washington} and
\cite[p.~261]{Lang}. If $(p, k)$ is an irregular pair, Washington examines the consequences of assuming that the probability that the cyclotomic unit $u_k$ is a $p$-th power is $1/p$, all independently of each other. Using this and heuristics on the distribution of $i_p$, he guesses that the number of counterexamples to Vandiver's conjecture for $p < x$ should grow like $\frac 12 \log\log x$. As an exercise in futility, we may try to refine Washington's heuristics using the known values of $i_p$, taking into account for example that the first irregular prime ($p = 37$) is unexpectedly large. We find that the expected number of counterexamples up to 12 million is about $0.674$, and that another $0.074$ counterexamples were expected between 12 million and 163 million (though of course we now know that there are none in either case). Many people believe that Vandiver's conjecture is true; it also seems reasonable
to believe that the conjecture is false but that the first counterexample
is so astronomically large that it may never be known. Similar remarks apply to the $\lambda_p = i_p$ conjecture.

As mentioned earlier, the table of irregular pairs, together
with check data for the Vandiver and cyclotomic invariant verifications,
can be found at \cite{URL}.

\section{Correctness}

Any computation on this scale is virtually guaranteed to encounter faults of various kinds, and this computation was no exception. On one occasion during the main Bernoulli number computation, a hardware or operating system failure --- we could not determine which --- resulted in several kilobytes of output being overwritten by random data, and it was necessary to repeat the computation for the affected primes. On another occasion, the Vandiver and cyclotomic runs located subtle memory errors (subsequently confirmed by other tests) on two workstations.  In comparing our data for those tests, we found a handful of examples (which were recalculated) were some sort of write error invalidating a very small number of lines. To combat this sort of occurrence, we took careful precautionary measures throughout the project, which we now describe.

The main irregular index computation was far too expensive to run more than once. We checked the correctness of this phase by several methods. During the main computation itself, for each $p$ we verified the identity 
\[
\sum_{k=0}^{p-3} 2^k (k + 1) B_k \equiv -4 \pmod p
\]
(see \cite[p.~720]{1-million}). Already this is a
strong consistency check, but it has the drawback that it leaves
no certificate that can be checked after the computation is finished. To increase our confidence in the output, we made our program record
somewhat more than just the irregular indices. Let $N_p
= \min(2\log p, \frac12(p-3))$. After computing all $B_k \bmod p$
for a given $p$, we extracted and stored the $N_p$ smallest pairs $(k, B_k \bmod
p)$, where the pairs are ordered first by $B_k \bmod p$ and then
by $k$. In particular, the irregular pairs are
listed first, and for the largest primes considered we store about
37 pairs. The resulting compressed file is 1.5 GB, and is available
for download \cite{URL}. We then used a completely independent
program and hardware to check these recorded pairs. For each $k$,
one may compute $B_k \bmod p$ in $O(p)$ operations using an identity such as \eqref{eq:voronoi}. We used the highly optimized implementation
described in \cite{multimodular}.

This scheme has an
additional property we might expect of a `certificate'. A skeptic
may randomly select and verify a pair $(k, B_k \bmod p)$ from our output file. Since presumably the only way to find a $k$ for which $B_k \bmod p$ is
`small' is to compute all the $B_k \bmod p$, the skeptic may be
persuaded that we must have actually performed the full computation
for each $p$. Note that recording \emph{all} the $B_k \bmod p$ that
we computed is impractical, requiring in the order of $10^6$ GB of
storage.

For the Vandiver tests, we computed $V_{p,k}$ twice for each irregular pair and cross-checked the results. This was done using programs written independently by the two authors, working separately and sharing no code. The programs were run on different hardware, using different libraries for the modular arithmetic (zn\_poly and NTL \cite{ntl}).

The first run used a new implementation of the method described in \cite{12-million}, updated for 64-bit hardware. For the second run, we used a slightly different algorithm. Instead of iterating over $c$ from $1$ up to
$(p-1)/2$, we write $c = 2^i g^j$, where $g$ is a primitive root
modulo $p$, and iterate over $0 \leq i < t$ in an inner loop and
$0 \leq j < (p-1)/t$ in an outer loop, where $t$ is the order
of $2$ modulo $p$. For each pair $(i, j)$ we only include the
multiplicand in the product if the proposed $c$ satisfies $1 \leq
c \leq (p-1)/2$. The advantage of this approach is that only a
single multiplication in $\ZZ/p\ZZ$ is required to update $c^{p-1-k}$
on each iteration, rather than the $O(\log p)$ operations needed
to compute each $c^{p-1-k}$ had the $c$ been processed sequentially.
Updating $z^c$ under this scheme requires only a single operation
in the inner loop (since $z^{2^{i+1} g^j} = (z^{2^i g^j})^2$), and
slightly more in the outer loop. The latter executes infrequently
for most $p$, since $t$ is usually large.

For the cyclotomic invariants, we again ran everything twice, using independent programs on different hardware. The first run used the implementation from \cite{12-million}, and the second was written from scratch using NTL's modular arithmetic.

Earlier computations of irregular pairs and associated data uncovered minor problems
in their predecessors' data,
so it was natural for us to compare our data with the results
up to 12 million from ten years ago.   
We found four errors in the online tables, and are confident that three of the cases
were introduced by an ill-advised editing, by hand, of the table to correct
cases in which a faulty computer had been used, and one of the cases
was output of the faulty computer itself.

\section{Hardware}

The aforementioned computations were performed on a number of different machines:

`Lonestar' is a 1300-node cluster at the Texas Advanced Computing Center (TACC), running a custom Linux operating system. Each node contains two dual-core 2.66 GHz Intel Xeon (Woodcrest) processors and has 8 GB RAM. This machine was used for computing irregular indices (about 119\,000 core-hours). We used various versions of GMP (4.2.1--4.2.4) \cite{gmp} for the large integer multiplication, patched with assembly code written by Jason Martin for the Core 2 architecture.

`Ranger' is a 3936-node cluster at TACC, running a custom Linux operating system. Each node contains four quad-core 2.3 GHz AMD Opteron (Barcelona) processors and has 32 GB RAM. At the time the computations were run, this machine was ranked by top500.org as the fourth most powerful supercomputer in the world. This machine was used for computing irregular indices (about 68\,000 core-hours). We used the same versions of GMP mentioned above, together with assembly code written by the second author and Torbj\"orn Granlund. This code has subsequently been incorporated into public releases of GMP; indeed, the needs of the present computation were a large part of the motivation for this work on GMP.

The first run of the Vandiver and cyclotomic invariant checks was performed on a network of 24 desktop machines, each with a 3.4 GHz Pentium 4 CPU (about 9\,000 core-hours). The second run was performed on a 16-core 2.6 GHz Opteron server (about 10\,000 core-hours). This machine was also used for the verification of the $B_k \bmod p$ (about 500 core-hours).

\section*{Acknowledgements}

Many thanks to Torbj\"orn Granlund for contributing his invaluable assembly programming expertise.

The authors acknowledge the Texas Advanced Computing Center (TACC) at The University of Texas at Austin for providing high-performance computing resources that have contributed to the research results reported within this paper.

Thanks to Fernando Rodriguez Villegas and the Department of Mathematics at the University of Texas at Austin for arranging access to the clusters at TACC, and to Salman Butt for his technical assistance with the clusters.

Thanks to the Department of Mathematics at Harvard University for supplying the Opteron server.

The first author thanks Richard Crandall and Amin Shokrollahi for useful conversations on this paper.

\bibliographystyle{amsalpha}
\bibliography{paper}

\newcommand{\etalchar}[1]{$^{#1}$}
\providecommand{\bysame}{\leavevmode\hbox to3em{\hrulefill}\thinspace}
\providecommand{\MR}{\relax\ifhmode\unskip\space\fi MR }
\providecommand{\MRhref}[2]{%
  \href{http://www.ams.org/mathscinet-getitem?mr=#1}{#2}
}
\providecommand{\href}[2]{#2}
\begin{thebibliography}{BCEM93}

\bibitem[BCE{\etalchar{+}}01]{12-million}
Joe Buhler, Richard Crandall, Reijo Ernvall, Tauno Mets{\"a}nkyl{\"a}, and
  M.~Amin Shokrollahi, \emph{Irregular primes and cyclotomic invariants to 12
  million}, J. Symbolic Comput. \textbf{31} (2001), no.~1-2, 89--96,
  Computational algebra and number theory (Milwaukee, WI, 1996).

\bibitem[BCEM93]{4-million}
J.~Buhler, R.~Crandall, R.~Ernvall, and T.~Mets{\"a}nkyl{\"a}, \emph{Irregular
  primes and cyclotomic invariants to four million}, Math. Comp. \textbf{61}
  (1993), no.~203, 151--153.

\bibitem[BCS92]{1-million}
J.~P. Buhler, R.~E. Crandall, and R.~W. Sompolski, \emph{Irregular primes to
  one million}, Math. Comp. \textbf{59} (1992), no.~200, 717--722.

\bibitem[BLS03]{tellegen}
Alin Bostan, Gr\'egoire Lecerf, and \'Eric Schost, \emph{Tellegen's principle
  into practice}, Symbolic and Algebraic Computation (J.~R. Sendra, ed.), ACM
  Press, 2003, Proceedings of ISSAC'03, Philadelphia, August 2003., pp.~37--44.

\bibitem[Blu70]{bluestein}
L.~Bluestein, \emph{A linear filtering approach to the computation of discrete
  {F}ourier transform}, Audio and Electroacoustics, IEEE Transactions on
  \textbf{18} (1970), no.~4, 451--455.

\bibitem[Cor08]{Corry}
Leo Corry, \emph{Fermat meets {SWAC}: {V}andiver, the {L}ehmers, computers, and
  number theory}, IEEE Ann. Hist. Comput. \textbf{30} (2008), no.~1, 38--49.

\bibitem[Dil09]{DilURL}
2009, http://www.mscs.dal.ca/\~{}dilcher/bernoulli.html/.

\bibitem[DSS91]{Dil}
Karl Dilcher, Ladislav Skula, and Ilja~Sh. Slavutski\v{\i}, \emph{Bernoulli
  numbers}, Queen's Papers in Pure and Applied Mathematics, vol.~87, Queen's
  University, Kingston, ON, 1991, Bibliography (1713--1990).

\bibitem[EM91]{EMa}
R.~Ernvall and T.~Mets{\"a}nkyl{\"a}, \emph{Cyclotomic invariants for primes
  between 125000 and 150000}, Math. Comp. \textbf{56} (1991), 851--858.

\bibitem[EM92]{EMb}
\bysame, \emph{Cyclotomic invariants for primes up to one million}, Math. Comp.
  \textbf{59} (1992), 249--250.

\bibitem[Gra08]{gmp}
Torbj\"orn Granlund, \emph{The {GNU} {M}ultiple {P}recision {A}rithmetic
  library}, 2008, http://gmplib.org/.

\bibitem[Har08a]{multimodular}
David Harvey, \emph{A multimodular algorithm for computing {B}ernoulli
  numbers}, to appear in Mathematics of Computation, preprint at
  http://arxiv.org/abs/0910.1926, 2008.

\bibitem[Har08b]{zn-poly}
David Harvey, \emph{The zn\_poly library}, 2008,
  http://www.cims.nyu.edu/\~{}harvey/zn\_poly/.

\bibitem[Har09a]{URL}
2009, http://www.cims.nyu.edu/\~{}harvey/irregular/.

\bibitem[Har09b]{cache-trunc-fft}
David Harvey, \emph{A cache-friendly truncated {FFT}}, Theoret. Comput. Sci.
  \textbf{410} (2009), no.~27-29, 2649--2658.

\bibitem[Har09c]{multipoint}
\bysame, \emph{Faster polynomial multiplication via multipoint {K}ronecker
  substitution}, J. Symbolic Comput. \textbf{44} (2009), no.~10, 1502--1510.

\bibitem[HQZ04]{middle-product}
Guillaume Hanrot, Michel Quercia, and Paul Zimmermann, \emph{The middle product
  algorithm, {I}.}, Appl. Algebra Engrg. Comm. Comput. \textbf{14} (2004),
  no.~6, 415--438.

\bibitem[Joh75]{John}
Wells Johnson, \emph{Irregular primes and cyclotomic invariants}, Math. Comp.
  \textbf{29} (1975), 113--120.

\bibitem[Lan90]{Lang}
Serge Lang, \emph{Cyclotomic fields {I} and {II}}, Graduate Texts in
  Mathematics, vol.~21, Springer-Verlag, New York, 1990.

\bibitem[Nus80]{nussbaumer}
Henri~J. Nussbaumer, \emph{Fast polynomial transform algorithms for digital
  convolution}, IEEE Trans. Acoust. Speech Signal Process. \textbf{28} (1980),
  no.~2, 205--215.

\bibitem[Sch77]{schonhage}
A.~Sch{\"o}nhage, \emph{Schnelle {M}ultiplikation von {P}olynomen \"uber
  {K}\"orpern der {C}harakteristik 2}, Acta Informat. \textbf{7} (1976/77),
  no.~4, 395--398.

\bibitem[Sho09]{ntl}
Victor Shoup, \emph{{NTL}: A library for doing number theory},
  http://www.shoup.net/ntl/, 2009.

\bibitem[Sla87]{Slav}
I.~Sh. Slavutski{\u\i}, \emph{A remark on the paper of {T}.\ {U}ehara: ``{O}n
  {$p$}-adic continuous functions determined by the {E}uler numbers'' [{R}ep.\
  {F}ac.\ {S}ci.\ {E}ngrg.\ {S}aga {U}niv.\ {M}ath.\ {N}o.\ 8 (1980), 1--8;
  {MR}0567622 (81e:12020)]}, Rep. Fac. Sci. Engrg. Saga Univ. Math. \textbf{15}
  (1987), 1--2.

\bibitem[SS71]{schonhage-strassen}
Arnold Sch{\"o}nhage and Volker Strassen, \emph{Schnelle {M}ultiplikation
  grosser {Z}ahlen}, Computing (Arch. Elektron. Rechnen) \textbf{7} (1971),
  281--292.

\bibitem[TW87]{WagTan}
Jonathan Tanner and Jr. Wagstaff, Samuel~S., \emph{New congruences for the
  {B}ernoulli numbers}, Math. Comp. \textbf{48} (1987), 341--350.

\bibitem[Van17]{Vnd}
H.~S. Vandiver, \emph{Symmetric functions formed by systems of elements of a
  finite algebra and their connection with fermat's quotient and {B}ernoulli
  numbers}, Ann. Math. \textbf{18} (1917), 105--114.

\bibitem[vdH04]{tft1}
Joris van~der Hoeven, \emph{The truncated {F}ourier transform and
  applications}, I{SSAC} 2004, ACM, New York, 2004, pp.~290--296.

\bibitem[vdH05]{tft2}
\bysame, \emph{Notes on the truncated {F}ourier transform}, unpublished,
  available from http://www.math.u-psud.fr/\~{}vdhoeven/, 2005.

\bibitem[Wag78]{Wag}
Jr. Wagstaff, Samuel~S., \emph{The irregular primes to 125000}, Math. Comp.
  \textbf{32} (1978), 583--591.

\bibitem[Was97]{Washington}
Lawrence~C. Washington, \emph{Introduction to cyclotomic fields}, second ed.,
  Graduate Texts in Mathematics, vol.~83, Springer-Verlag, New York, 1997.

\end{thebibliography}

\end{document}